\documentclass[11pt,a4paper]{article}
\usepackage[T1]{fontenc}
\usepackage{lmodern,microtype}
\usepackage[hmargin=26mm,vmargin=25mm]{geometry}
\usepackage{amsmath,amssymb,amsthm,mathtools}
\usepackage[hidelinks]{hyperref}
\newtheorem{theorem}{Theorem}[section]
\newtheorem{proposition}[theorem]{Proposition}
\newtheorem{corollary}[theorem]{Corollary}

\numberwithin{equation}{section}
\newcommand{\Es}{E_{\rm sym}}
\newcommand{\F}{\mathcal F_\omega}
\newcommand{\T}{\mathbb T}
\newcommand{\R}{\mathbb R}
\newcommand{\dv}{\,dV_g}

\newcommand{\HS}{\mathrm{HS}}
\newcommand{\rank}{\operatorname{rank}}
\newcommand{\osc}{\operatorname{osc}}
\newcommand{\Harm}{\mathcal H_g}
\DeclareMathOperator{\tr}{tr}
\title{Rank-sensitive Hodge bounds for symphonic energy\\and minimizers on sheared tori}
\author{Xiangzhi Cao\thanks{aaa7756kijlp@163.com}\,\,\thanks{Nanjing Xiaozhuang University}}
\date{}
\begin{document}
\maketitle
\begin{abstract}
We ask how far the minimum of symphonic energy in a homotopy class can exceed the Hodge lower bound from a pulled-back two-form. For projections from a family of sheared four-tori to the flat two-torus, we compute the minimum and the gap for every smooth shear profile. The gap depends on the variance of the profile. Attainment of the projection's pullback area energy by a coclosed test depends on its oscillation. We determine exactly which variances allow two profiles with the same positive gap and different outcomes for this test, and compute the convex relaxation beyond the oscillation threshold. In arbitrary source dimension the optimized test improves the Hodge bound precisely when the harmonic representative exceeds the rank allowed for the pullback on a set of positive measure. We also prove low-degree cohomology estimates, a sharp universal degree exponent for targets detected by low-degree cup products, and an exact defect identity that forces an immersed source to be symplectic at equality.
\end{abstract}

\section{Cohomology and the energy gap}
Let $M$ be a closed Riemannian manifold, not necessarily orientable, and $N$ a Riemannian manifold, not necessarily compact. All integrals on $M$ use its Riemannian density. For a smooth map $u:M\to N$ and a smooth closed two-form $\omega$ on $N$, write
\[
 \beta=u^*\omega,\qquad \eta=\Harm[u^*\omega],\qquad
 \F(u)=\tfrac12\|\beta\|_2^2,\qquad
 \Es(u)=\tfrac14\int_M\tr((du\,du^*)^2)\dv .
\]
Here $\Harm$ selects the harmonic representative in real de Rham cohomology and $du^*$ is the metric adjoint. Critical points of $\Es$ are symphonic maps \cite{NT}. The pullback class is fixed in a homotopy class. How close can its harmonic representative bring us to the actual minimum of the energy?

Consider first the projection $p=(x,y):\T^4\to\T^2$ to the unit flat torus with area form $\omega$. On the source take the metric for which $dx,dy,dz+b\,dy,dt+f(z)\,dx$ are orthonormal. Here $b$ is real, $f$ is smooth of period one, and $V$ is its variance on one period. Theorems~\ref{thm:torus} and~\ref{thm:certificate} give
\[
 \min_{u\simeq p}\Es(u)=\tfrac12,\qquad
 \min_{u\simeq p}\Es(u)-\tfrac12\|\Harm[p^*\omega]\|_2^2
       =\frac{b^2V}{2(1+b^2)(1+V)}.
\]
For the area functional $\F$, a coclosed test realizes $\F(p)=1/2$ exactly when $|b|(\max f-\min f)\le2$. Consequently the same symphonic minimum and the same positive gap occur on opposite sides of this threshold. For fixed $b\ne0$ this happens precisely for $0<V<b^{-2}$ (Corollary~\ref{cor:variance}). Beyond the threshold Proposition~\ref{cor:strictdual} computes the best test bound. The infimum of $\F$ over maps in that regime remains open.

The mechanism extends to arbitrary dimensions. Hodge orthogonality gives
\begin{equation}\label{eq:orth}
 2\F(u)=\|\eta\|_2^2+\|\beta-\eta\|_2^2.
\end{equation}
Speight and Svensson \cite{SS} obtain homotopy minimality for $\F$ from a coclosed pullback. A surface area form, for example, has pullbacks of rank at most two, while $\eta$ can have higher rank. Theorem~\ref{thm:test} turns such a rank violation into a uniform positive gap using coclosed test forms. The pointwise norm is a skew-symmetric instance of the dual spectral $k$-support norm \cite{MPS}. Proposition~\ref{prop:quotient} identifies the optimal bound as a convex minimum over forms. Slobodeanu and Svensson \cite{SlSv} studied the second variation at surface-valued submersions with nonintegrable horizontal distributions. Hattori \cite{Hat} established lower bounds for maps by calibration, including examples on flat tori. Here the sheared source metric allows direct computation of both the harmonic representative and the optimal test bound.

Han and Wei \cite{HW} study the same functional under the name $\Phi$-energy: their density is $\tfrac14\sum_{i,j}h(du(e_i),du(e_j))^2$ in a source orthonormal frame. For compact manifolds their superstrong instability hypothesis on either endpoint yields arbitrarily small energy in each homotopy class. They also prove that a compact $\Phi$-superstrongly unstable manifold is four-connected. Proposition~\ref{prop:low} gives a complementary obstruction applicable to individual classes: a nonzero real pullback in any degree from one to four forces a positive energy infimum. Its proof uses exterior powers and H\"older's inequality. Van Schaftingen \cite{VS} gives related estimates for induced cohomology maps using critical Sobolev energies and gap potentials. Cup products then give the degree bound in Corollary~\ref{cor:degree}. For almost Hermitian targets, Theorem~\ref{thm:defects} separates the harmonic and pointwise algebraic defects and describes the geometry at equality.

\section{Cohomology actions and mapping degree}\label{sec:low}
Let $M^m,N^n$ be closed Riemannian manifolds with fixed metrics. Equip real cohomology with the $L^2$ inner product on harmonic representatives, using the Riemannian densities on both manifolds. Orientations are needed below only for the signed mapping degree. Write $e=\frac14\tr((du\,du^*)^2)$.

\begin{proposition}[Low-degree cohomology]\label{prop:low}
For $1\le k\le\min\{4,m,n\}$ there is a constant $C_k$ depending only on the two fixed metrics and manifolds such that
\begin{equation}\label{eq:low}
 \big\|u^*:H^k(N;\R)\longrightarrow H^k(M;\R)\big\|_{\HS}
 \le C_k\Es(u)^{k/4}.
\end{equation}
For any finite energy bound the induced homomorphisms on integral cohomology in these degrees form a finite set, both modulo torsion and on the full groups.
\end{proposition}
\begin{proof}
Take orthonormal harmonic bases $\{\alpha_i\}$ on $N$ and $\{\theta_j\}$ on $M$. The matrix entries are $F_{ji}=\int_M\langle u^*\alpha_i,\theta_j\rangle\dv$, since each $\theta_j$ is orthogonal to the exact part of a closed form. If $s_1\ge\cdots\ge s_{\min(m,n)}$ are the singular values of $du$, the norm of its $k$th exterior power satisfies
\[
 s_1\cdots s_k\le\left(k^{-1}\sum_{\ell=1}^k s_\ell^4\right)^{k/4}
 \le (4e/k)^{k/4}.
\]
H\"older's inequality, including the endpoint $k=4$, gives
\[
 |F_{ji}|\le (4/k)^{k/4}\|\alpha_i\|_\infty\|\theta_j\|_\infty
 \operatorname{Vol}_g(M)^{1-k/4}\Es(u)^{k/4}.
\]
Taking the Euclidean norm of these finitely many entries proves \eqref{eq:low}. The images of integral cohomology modulo torsion are full lattices in the two real spaces. Homomorphisms between these lattices are integer matrices in integral bases, and every bounded subset contains finitely many of them. Two integral homomorphisms inducing the same real map differ by a homomorphism from $H^k(N;\mathbb Z)$ to the torsion subgroup of $H^k(M;\mathbb Z)$. The domain is finitely generated and the latter subgroup is finite, so there are only finitely many such differences. The energy bound is simultaneous in the indicated degrees.
\end{proof}

\noindent In particular, if a fixed homotopy class contains maps whose symphonic energies tend to zero, its real pullback vanishes in every degree covered by \eqref{eq:low}. A nonzero pullback in one such degree gives the quantitative obstruction $\inf_{v\simeq u}\Es(v)\ge(\|u^*\|_{\HS}/C_k)^{4/k}$ for $C_k>0$. This obstruction uses real cohomology, even though the preceding finiteness statement applies to integral cohomology.

\begin{corollary}[Degree detected by cup products]\label{cor:degree}
Suppose $M,N$ are connected, closed and oriented with $m=n\ge1$. Assume there are classes $a_j\in H^{k_j}(N;\R)$ with $1\le k_j\le4$, $\sum_j k_j=m$, and $I_N=\langle\smile_j a_j,[N]\rangle\ne0$. There is a constant $c>0$ depending only on the chosen data such that
\begin{equation}\label{eq:cup}
 \Es(u)\ge c|\deg u|^{4/m}.
\end{equation}
The hypothesis holds for tori of every dimension and for closed symplectic targets of every even dimension.
\end{corollary}
\begin{proof}
By naturality $(\deg u)I_N=\langle\smile_j u^*a_j,[M]\rangle$. The multilinear cup pairing on the finite-dimensional Hodge spaces of $M$ has finite norm $K_M$. Applying \eqref{eq:low} to each factor yields
\[
 |\deg u|\,|I_N|\le K_M\prod_j(C_{k_j}\|a_j\|)\,\Es(u)^{m/4}.
\]
If the coefficient on the right vanishes, all maps have degree zero. Otherwise rearrangement gives \eqref{eq:cup}. On a torus take the product of the coordinate one-classes. On a closed symplectic $2d$-manifold the $d$th power of the symplectic class pairs nontrivially with the fundamental class.
\end{proof}

On the unit-volume square torus $\T^m$, the multiplication maps $x\mapsto jx$ have degree $j^m$ and energy $mj^4/4$ for positive integers $j$. Thus $4/m$ is the largest possible exponent in an estimate of this form valid throughout the stated class. For comparison, when $m>4$, rescale a smooth degree-one map $B^m/\partial B^m\to S^m$, constant on a collar, to a coordinate ball of radius $\varepsilon$ in $S^m$ and extend it by that constant. Its degree stays one while its energy is $O(\varepsilon^{m-4})$. This shows why some restriction on the high-dimensional target matters.

\section{The energy defect for an almost Hermitian target}\label{sec:defects}
Return to a closed source $M$ of arbitrary dimension. On $N$ let $J$ be an orthogonal almost complex structure with fundamental two-form $\omega(v,w)=h(Jv,w)$. For a single map we require only $d(u^*\omega)=0$. This holds for every map to an almost K\"ahler target, and for any map of pointwise rank at most two since $d(u^*\omega)=u^*(d\omega)$. Set $A=du$ and $B=AA^*$ on the pulled-back target tangent bundle. All endomorphism norms below are Hilbert--Schmidt norms.

\begin{theorem}[Two squared defects]\label{thm:defects}
For $\eta=\Harm[u^*\omega]$ one has
\begin{equation}\label{eq:defects}
 \Es(u)=\tfrac12\|\eta\|_2^2+
 \tfrac12\|u^*\omega-\eta\|_2^2+\tfrac18\|BJ-JB\|_2^2 .
\end{equation}
Equality in $\Es(u)\ge\|\eta\|_2^2/2$ holds precisely when $u^*\omega$ is harmonic and $B$ commutes with $J$ everywhere. An equality map minimizes energy among smooth maps with closed pullback in the same cohomology class. If $d\omega=0$, it minimizes in its homotopy class and is a stable symphonic map.
\end{theorem}
\begin{proof}
Since $J^*=-J$ and $J^2=-I$, the squared form norm and the squared commutator norm are
\[
 |u^*\omega|^2=\tfrac12\|A^*JA\|_{\HS}^2=-\tfrac12\tr(BJBJ),\qquad
 \|BJ-JB\|_{\HS}^2=2\tr(B^2)+2\tr(BJBJ).
\]
The second identity follows by squaring the self-adjoint commutator and cycling the trace. Eliminating $\tr(BJBJ)$ gives the pointwise density $\frac12|u^*\omega|^2+\frac18\|BJ-JB\|_{\HS}^2$. Integrating and applying \eqref{eq:orth} to the closed pullback proves \eqref{eq:defects}. Both remainders are nonnegative, and the harmonic term is fixed throughout the specified class. If $d\omega=0$, homotopy invariance preserves that class along every smooth variation. An equality map then has vanishing first variation and nonnegative second variation.
\end{proof}

\begin{corollary}[Topology at equality]\label{cor:horizontal}
For a map attaining the bound in Theorem~\ref{thm:defects}, $\rank du$ is even at each point. If the map has constant rank, its horizontal bundle $H=(\ker du)^\perp$ admits a complex vector bundle structure. If it is an immersion, then $u^*\omega$ is a symplectic form on $M$.
\end{corollary}
\begin{proof}
At equality $J$ preserves $\operatorname{im}B=\operatorname{im}A$, so its real dimension is even. At constant rank this image is a smooth subbundle. Transport $J$ to $H$ by the bundle isomorphism $du|_H$. For an immersion $H=TM$. If $v\ne0$ and $w=(du|_H)^{-1}Jdu(v)$, then $(u^*\omega)(v,w)=|du(v)|^2>0$. Thus the closed pullback form is nondegenerate.
\end{proof}
The commutation condition is the pseudo horizontally weakly conformal condition studied by Loubeau \cite{Lou} and Loubeau and Mo \cite{LM}. The complex structure on the horizontal bundle is the restriction of the associated $f$-structure of \cite{Lou}. The energy identity specifies when equality enforces this structure. In the sheared projection of Theorem~\ref{thm:torus}, $dp\,dp^*=I_2$, so \eqref{eq:defects} identifies the entire excess in \eqref{eq:min} with the squared distance of the pulled-back area form from its harmonic representative, divided by two.

\section{Rank bounds in arbitrary dimension}
Let $m=\dim M$, $q=\lfloor m/2\rfloor$. A two-form $\theta$ has an orthonormal normal form $\sum_{i=1}^q a_i e^{2i-1}\wedge e^{2i}$, where $a_1\ge\cdots\ge a_q\ge0$. We use $|\theta|^2=\sum a_i^2$ and put
\[
 t_r(\theta)^2=\sum_{i\le\min(r,q)}a_i^2,\qquad
 T_r(\theta)=\int_M t_r(\theta)^2\dv,\qquad S(\theta)=\|\theta\|_2^2.
\]
For $r=1$, $t_1$ is the comass. If there are no two-forms or the pullback rank vanishes, the estimates concerning nonzero classes are void. In what follows $m\ge2$ and $r\ge1$.

\begin{theorem}\label{thm:test}
Suppose $\rank(u^*\omega)\le 2r$ almost everywhere. For every smooth coclosed two-form $\psi$ on $M$ with $T_r(\psi)>0$,
\begin{equation}\label{eq:test}
 2\F(u)\ \ge\ \frac{|\langle\eta,\psi\rangle_{L^2}|^2}{T_r(\psi)} .
\end{equation}
In particular, if $S=S(\eta)>0$ and $T=T_r(\eta)$, then $0<T\le S$ and
\begin{equation}\label{eq:gain}
 2\F(u)\ \ge\ \frac{S^2}{T}
 = S+(S-T)+\frac{(S-T)^2}{T}.
\end{equation}
If $\omega$ has comass at most $K>0$ along the image of $u$, then the same right-hand sides bound $2K^2\Es(u)$ from below. If $\omega$ has rank at most $2r$ and comass at most $K$ on $N$, these estimates hold uniformly over all smooth maps with the same class $[u^*\omega]$.
\end{theorem}
\begin{proof}
At a point, write $\beta=\sum_{i=1}^s b_i v_i^\flat\wedge w_i^\flat$ in orthonormal pairs with $s\le r$. If $\Psi$ is the skew operator of $\psi$, then
\[
 \sum_{i=1}^s|\psi(v_i,w_i)|^2
 \le \tfrac12\sum_{i=1}^s(|\Psi v_i|^2+|\Psi w_i|^2)
 \le \sum_{i\le \min(r,q)}a_i(\psi)^2.
\]
For the last step diagonalize $\Psi^*\Psi$: the trace on a $2s$-plane is no larger than the sum of its largest $2s$ eigenvalues, with each block value $a_i^2$ counted twice. Cauchy--Schwarz gives $|\langle\beta,\psi\rangle|\le|\beta|t_r(\psi)$. This coefficient is optimal pointwise: the truncation of $\psi$ to its first $r$ blocks realizes equality. Since $\beta-\eta$ is exact, it pairs to zero with $\psi$ on the closed source. An integral Cauchy--Schwarz inequality proves \eqref{eq:test}. Take $\psi=\eta$ and use that nonzero $\eta$ has positive largest block on a set of positive measure. This proves \eqref{eq:gain}.

For the metric comparison, write $A=du$, $B=AA^*$, and $\omega(v,w)=h(\Omega v,w)$. In an orthonormal eigenbasis of $B$ with eigenvalues $b_i\ge0$,
\[
 2|\beta|^2=\sum_{i,j}b_i b_j\Omega_{ij}^2
 \le\sum_i b_i^2\sum_j\Omega_{ij}^2
 \le K^2\tr(B^2).
\]
The first inequality uses $2b_i b_j\le b_i^2+b_j^2$ and skew symmetry; the second uses the operator norm bound on $\Omega$. Integration yields $\F(u)\le K^2\Es(u)$. The class and harmonic representative are unchanged across the stated family.
\end{proof}

\begin{corollary}\label{cor:rank}
For a fixed nonzero realizable pullback class, the harmonic-test improvement in \eqref{eq:gain} is strictly positive if and only if $\rank\eta>2r$ on a set of positive measure. In that case the excess of $\F$ above $S/2$ is at least $S(S-T)/(2T)$ for every map in the class subject to the rank cap.
\end{corollary}
\begin{proof}
$S-T$ is the integral of the sum of squares of the discarded blocks. Subtract $S$ from $S^2/T$. For a smooth harmonic form a rank violation at one point persists on an open set.
\end{proof}

No almost complex structure, positivity of curvature, nonzero degree, constant rank or lower bound on the differential appears here. In fact \eqref{eq:test} itself is a Hilbert-space statement for an $L^2$ form $\beta$ of bounded rank: the only cohomological input is $\beta-\eta$ in the $L^2$-closure of exact forms. It consequently applies without asserting any density theorem for manifold-valued Sobolev maps. One can use a boundary provided the orthogonality to the selected tests is imposed. Without a boundary condition integration by parts has a boundary term. No such extension is needed for the examples.

For comparison, the squared distance of $\eta(x)$ to the rank-$2r$ cone is $|\eta(x)|^2-t_r(\eta(x))^2$. Indeed the pointwise support estimate applied to any admissible $\xi$ gives
$|\eta-\xi|^2\ge|\eta|^2-t_r(\eta)^2+(|\xi|-t_r(\eta))^2$, and truncation attains the minimum. Using this distance alone in \eqref{eq:orth} gives $2\F\ge S+(S-T)$ and loses the last term in \eqref{eq:gain}. Neither bound is asserted to be the exact minimum over pullbacks. Equality in the harmonic test for $S>0$ requires, and is equivalent to, pointwise support equality with nonnegative pairing and $|\beta|=(S/T)t_r(\eta)$ almost everywhere. This includes possible ties among the block values.

Define the optimized lower bound
\begin{equation}\label{eq:dual}
 D_r(c)=\sup_{\substack{\delta\psi=0\\T_r(\psi)>0}}
 \frac{|\langle\Harm c,\psi\rangle|^2}{T_r(\psi)}.
\end{equation}
Smooth tests suffice by definition. Since $t_r$ dominates a positive dimensional multiple of the Euclidean norm, the supremum is finite. The formula provides a certificate whenever one test attains $2\F(u)$; it is not a general duality assertion for the nonconvex minimization over maps.

There is, however, an exact convex problem behind this certificate. Let $t_r^\circ$ be the pointwise polar norm:
$t_r^\circ(\xi)=\sup_{t_r(\psi)\le1}\langle\xi,\psi\rangle$.
Its unit ball is the convex hull of the unit Euclidean ball restricted to forms of rank at most $2r$, by the support calculation in the proof above and finite-dimensional polarity. In particular $t_r^\circ(\xi)=|\xi|$ on that rank cone. Write $X$ for the $L^2$-closure of smooth exact two-forms.

\begin{proposition}[Exact dual problem]\label{prop:quotient}
For every real two-class $c$, with $\eta=\Harm c$, one has
\[
 D_r(c)=\min_{\gamma\in\eta+X}\int_M t_r^\circ(\gamma)^2\dv .
\]
If $S(\eta)>0$, then $D_r(c)>S(\eta)$ if and only if the harmonic representative has rank greater than $2r$ on a set of positive measure.
\end{proposition}
\begin{proof}
Give $L^2(\Lambda^2M)$ the equivalent norm $n(\gamma)=\|t_r^\circ(\gamma)\|_2$.
Its dual norm, in the ordinary integrated pairing, is $\sqrt{T_r(\psi)}$: the upper bound is polar Cauchy--Schwarz, and the reverse follows by choosing at each point a unit polar-norm maximizer and scaling it by $t_r(\psi)$. One can use measurable approximate maximizers instead: in measurable orthonormal local frames take a countable dense set in the polar unit ball and select the first candidate within a fixed relative error of the pointwise supremum. The resulting fields, after scaling, are in $L^2$ by norm equivalence. Letting the error decrease proves the norm identity.
The quotient norm of $\eta+X$ is the supremum of its pairings with the unit dual ball in $X^\perp$, by Hahn--Banach. Here $X^\perp$ consists of distributionally coclosed forms. Smooth coclosed tests are dense there: Hodge heat regularization commutes with the codifferential and converges strongly in $L^2$, hence also in the dual norm. Squaring gives the claimed equality. The minimum exists by weak compactness of bounded sets in $L^2$ and weak lower semicontinuity of the norm on the closed affine subspace.
For rank excess use \eqref{eq:gain}. Without rank excess $\eta$ itself is an admissible competitor of cost $S(\eta)$, and the harmonic test gives the matching lower bound.
\end{proof}
The proposition identifies a relaxation over forms, not its realizability by differentials of maps. For $r=1$ its polar norm is the mass norm dual to comass; see Bangert and Katz \cite{BK} for related cohomological norms in systolic geometry.

Under the hypotheses of Theorem~\ref{thm:defects}, combining the pointwise identity in its proof with \eqref{eq:test} gives, whenever $\rank\beta\le 2r$,
\[
 \Es(u)\ge \tfrac12D_r([\beta])+\tfrac18\|BJ-JB\|_2^2 .
\]
Thus the improved cohomological bound retains the full commutator contribution to the energy.

\section{An exact minimum with nonharmonic pullback}
Let $\T^j=\R^j/\mathbb Z^j$ and write $x,y,z,t$ for unit-period coordinates. Given $b\in\R$ and a smooth periodic function $f=f(z)$, define an orthonormal coframe for $g_{b,f}$ by
\begin{equation}\label{eq:coframe}
 e^1=dx,\quad e^2=dy,\quad e^3=dz+b\,dy,\quad e^4=dt+f(z)\,dx .
\end{equation}
Use the orientation $dx\wedge dy\wedge dz\wedge dt$. The volume density is independent of $b,f$ and has total mass one. Give $\T^2$ the unit square metric and area form $\omega$. The projection $p=(x,y)$ has $\beta_0=p^*\omega=e^{12}$. Put
\[
 a=\int_0^1 f\,dz,\qquad V=\int_0^1(f-a)^2\,dz,\qquad L=1+b^2 .
\]
The notation $e^{ij}$ means $e^i\wedge e^j$.

\begin{theorem}\label{thm:torus}
For all $b,f$, the projection is, up to translations of the target, the unique smooth minimizer of $\Es$ in its homotopy class. In particular it is a stable symphonic map. The exact excess over the Hodge bound is
\begin{equation}\label{eq:min}
 \min_{u\simeq p}\Es(u)=\tfrac12,\qquad
 \min_{u\simeq p}\Es(u)-\tfrac12\|\Harm[\beta_0]\|_2^2
 =\frac{b^2 V}{2(1+b^2)(1+V)} .
\end{equation}
The harmonic representative has rank four on a nonempty open set precisely when $b\ne0$ and $f$ is nonconstant. In this family that condition is also necessary and sufficient for a positive gap above the Hodge bound.
\end{theorem}

We first record a quantitative version of the convexity argument. If $A$ is a real matrix, put $P(A)=\frac14\tr((AA^*)^2)$.
\begin{proposition}\label{prop:convex}
If $v:M\to\T^k$ is a smooth symphonic map to a flat torus and $u\simeq v$, let $\varphi:M\to\R^k$ be a lift of the nullhomotopic difference $u-v$ in flat target coordinates. Then
\begin{equation}\label{eq:coercivity}
 \Es(u)-\Es(v)\ge \tfrac13\int_M P(d\varphi)\dv.
\end{equation}
In particular any two smooth symphonic maps in a common homotopy class differ by a translation when $M$ is connected.
\end{proposition}
\begin{proof}
In orthonormal frames $DP(A)[H]=\langle AA^*A,H\rangle_{\HS}$, and
\[
 D^2P(A)[H,H]=\tfrac12|A^*H+H^*A|_{\HS}^2+\tr(A^*A H^*H)\ge0 .
\]
For $H\ne0$, $d^2P(A+sH)/ds^2$ is a nonnegative quadratic polynomial on $\R$ with leading coefficient $12P(H)$. It dominates $12P(H)(s-s_0)^2$ for its vertex $s_0$. Since $\min_z\int_0^1(1-s)(s-z)^2ds=1/36$, Taylor's integral formula gives
$P(A+H)-P(A)-DP(A)[H]\ge P(H)/3$. The constant is sharp for matrices: take $A=-H/3$. Integrate with $A=dv,H=d\varphi$; stationarity of $v$ annihilates the linear term by the variation $v+s\varphi$ modulo the target lattice. The last assertion follows by interchanging $v$ and $u$ in the inequality.
\end{proof}

\begin{proof}[Proof of Theorem~\ref{thm:torus}]
In the dual frame the rows of $dp$ are $(1,0,0,0)$ and $(0,1,0,0)$, so $dp\,dp^*=I_2$ and $P(dp)=1/2$. The linear term in the first variation integrates to
\[
 \int_{\T^4}((\partial_x-f\partial_t)\varphi_1+
                    (\partial_y-b\partial_z)\varphi_2)\dv=0
\]
by periodicity. The two displayed vector fields are divergence-free, including for nonconstant $f(z)$. Proposition~\ref{prop:convex} now gives the minimum and uniqueness.

To compute $\eta=\Harm[\beta_0]$, average over the isometric translations in $x,y,t$. They act trivially on cohomology; uniqueness makes $\eta$ invariant under all of them. An invariant exact two-form has an invariant primitive by averaging. Consequently the invariant closed representatives in this class are exactly
\begin{equation}\label{eq:rep}
 \gamma=dx\wedge dy+
 dz\wedge(P_0(z)dx+Q_0(z)dy+R(z)dt),
\end{equation}
where all three coefficients have zero mean. Conversely zero-mean coefficients have periodic primitives. Set $U=P_0-fR$. In the orthonormal frame this gives
\[
 \gamma=(1+bU)e^{12}-Ue^{13}-Q_0e^{23}+R e^{34}-bR e^{24}.
\]
Thus its squared norm integrates to the integral of
$(1+bU)^2+U^2+Q_0^2+L R^2$. We minimize this positive definite quadratic functional subject to $\int(U+fR)=\int R=\int Q_0=0$. At its minimizer $Q_0=0$; variation in $P_0$ gives $b+LU$ constant, and variation in $R$ at fixed $P_0$ gives $LR-f(b+LU)$ constant. The mean constraints yield
\begin{equation}\label{eq:solution}
 U=-\frac{bV}{L(1+V)},\qquad R=\frac{b(f-a)}{L(1+V)},\qquad P_0=U+fR .
\end{equation}
These smooth coefficients meet the constraints and solve the variational equations, hence give the harmonic representative. By \eqref{eq:orth}, $\|\eta\|_2^2=\langle\eta,\beta_0\rangle_{L^2}=1+bU$. This proves \eqref{eq:min}. In coordinates $\eta\wedge\eta=2R\,dx\wedge dy\wedge dz\wedge dt$, proving the rank criterion.
\end{proof}

\section{When a rank-one test gives the exact area energy}
Here ``rank-one test'' refers to the support function $t_1$, not to the rank of the test form: it may have rank four. Keep the data of \eqref{eq:coframe}.
\begin{theorem}\label{thm:certificate}
There is a smooth coclosed test $\psi$ for which \eqref{eq:test} certifies $2\F(p)=1$ if and only if
\begin{equation}\label{eq:osc}
 |b|\,\osc(f)\le2,\qquad \osc(f)=\max f-\min f.
\end{equation}
In this range $\min\F=1/2$ over all maps with $[u^*\omega]=[\beta_0]$, and also over the homotopy class of $p$. Outside this range no single smooth coclosed test attains the value for $p$. The symphonic minimum in \eqref{eq:min} holds on both sides of the threshold.
\end{theorem}
\begin{proof}
A test attaining equality in the two inequalities used for \eqref{eq:test} must have constant comass since $|\beta_0|=1$. After scalar normalization and a sign choice it has comass one and $\langle e^{12},\psi\rangle=1$ everywhere. The skew operator of $\psi$ then maps the oriented $12$-plane to itself: the $12$ coefficient saturates its operator norm, and no mixed coefficient in those two columns is possible. It follows, and conversely suffices for pointwise saturation, that
\[
 \psi=e^{12}+\lambda e^{34},\qquad |\lambda|\le1.
\]
Coclosedness is equivalent to $d(*\psi)=d(e^{34}+\lambda\,dx\wedge dy)=0$. Expanding $e^{34}$ using \eqref{eq:coframe} leaves precisely
\[
 d\big((\lambda-bf)dx\wedge dy\big)=0.
\]
Thus $\lambda=bf(z)+h(x,y)$ for a periodic function $h$. Such a $\lambda$ has absolute value at most one for a suitable $h$ if and only if the interval of values of $bf$ has length at most two. Its negative midpoint gives a constant $h$ in the admissible range. The resulting test gives $\langle\eta,\psi\rangle=\langle\beta_0,\psi\rangle=1$ and $T_1(\psi)=1$. The projection realizes the bound. The reverse implication concerned all smooth tests, not merely translation-invariant ones.
\end{proof}

\begin{proposition}[The value on both sides of the threshold]\label{cor:strictdual}
Use the notation $\langle h\rangle=\int_0^1h(z)\,dz$ and $\sigma(t)=\operatorname{sign}(t)(|t|-1)_+$, with $\sigma(0)=0$. There are real numbers $\alpha,\mu$ such that
\[
 R_*(z)=L^{-1}\sigma(\alpha f(z)-\mu),\quad
 \langle R_*\rangle=0,\quad
 \alpha=b-Lw,\qquad w=\langle fR_*\rangle .
\]
The resulting $R_*$ is unique, and the optimized bound is
\begin{equation}\label{eq:clipped}
 D_1([\beta_0])=1-L\big(w^2+\langle R_*^2\rangle\big).
\end{equation}
In particular it is strictly less than one when $|b|\osc(f)>2$.
\end{proposition}
\begin{proof}
Use the convex problem in Proposition~\ref{prop:quotient}. Averaging over $x,y,t$ preserves the affine class and does not increase the integral, by convexity and translation invariance. The invariant competitors are \eqref{eq:rep} with zero-mean $L^2$ coefficients: average smooth exact approximants for one inclusion, and approximate each zero-mean coefficient by smooth zero-mean functions for the other. For a two-form in four dimensions the polar of comass is the sum of its two block values. Indeed a test of comass one bounds each normal block separately, and the sum of the two oriented unit blocks attains this bound. Its square on $\gamma$ is therefore $|\gamma|^2+2|R|$, since $\gamma\wedge\gamma=2R\,dV_g$. For fixed $R$, the optimal $Q_0$ is zero and the optimal $U$ is its prescribed mean $-w$ by strict convexity. The remaining cost is
\[
 J(R)=1-2bw+Lw^2+L\langle R^2\rangle+2\langle|R|\rangle ,
 \qquad w=\langle fR\rangle,\quad \langle R\rangle=0 .
\]
This coercive strictly convex functional has a unique minimizer in the closed mean-zero subspace of $L^2$. Its subgradient equation is
$LR_*+\xi=\alpha f-\mu$, where $\xi\in\partial|R_*|$ pointwise. Here $\mu$ is constant since the orthogonal complement of the mean-zero subspace consists of constants; the sum rule applies since all terms are continuous on $L^2$ of the unit circle. Inverting the scalar inclusion gives the displayed formula for $R_*$. Multiplying by $R_*$ and averaging yields $L\langle R_*^2\rangle+\langle |R_*|\rangle=(b-Lw)w$, hence \eqref{eq:clipped}. The zero solution is admissible in the inclusion exactly when $|bf-\mu|\le1$ everywhere for some constant $\mu$, equivalently \eqref{eq:osc}. Otherwise $\langle R_*^2\rangle>0$.
\end{proof}

\begin{corollary}[Shear profiles of equal variance]\label{cor:variance}
Fix $b\ne0$ and $V\ge0$. There exist two smooth mean-zero profiles of variance $V$ with the same symphonic minimum and Hodge bound, exactly one of which admits a coclosed test realizing $\F(p)$, if and only if $0<V<b^{-2}$. For each such pair the values of $D_1([\beta_0])$ are different.
\end{corollary}
\begin{proof}
Write $\ell=\min f$, $h=\max f$ and $a=\langle f\rangle$. A continuous nonconstant profile takes values strictly between its extrema on an open set. Integrating $(f-\ell)(h-f)$ therefore gives $V<(h-a)(a-\ell)\le(h-\ell)^2/4$. A profile of positive variance satisfying the test condition must have $V<b^{-2}$. At $V=0$ every profile is constant and satisfies it.

Conversely suppose $0<V<b^{-2}$. Set $q_t(z)=\tanh(t\sin(2\pi z))$. Translation by a half-period changes its sign, so its mean vanishes. By dominated convergence $\langle q_t^2\rangle\to1$ as $t\to\infty$. The profile $\sqrt{V}\,q_t/\|q_t\|_2$ has variance $V$ and oscillation tending to $2\sqrt{V}<2/|b|$, so it admits a saturating test for large $t$.

To obtain the other profile, choose a smooth bump $\chi$ supported in $(0,1)$, with minimum zero and maximum one. On the unit circle put $g_\varepsilon(z)=\chi(z/\varepsilon)$ on $(0,\varepsilon)$ and zero elsewhere. Subtract its mean and multiply by $\sqrt{V}/\|g_\varepsilon-\langle g_\varepsilon\rangle\|_2$. This produces variance $V$. Its oscillation tends to infinity because the unscaled oscillation is one and the denominator is at most $\sqrt{\varepsilon}\|\chi\|_2$. For small $\varepsilon$ the threshold is exceeded. Formula~\eqref{eq:min} gives the same minimum and harmonic norm in the two metrics. By Theorem~\ref{thm:certificate} and Proposition~\ref{cor:strictdual}, $D_1=1$ for the first metric and $D_1<1$ for the second.
\end{proof}

For $b=1$ and $f(z)=s\cos(2\pi z)$, the threshold is $|s|=1$. At $s=1$ both minima are $1/2$, whereas the Hodge bound is $5/12$. For $s=2$ the symphonic minimum is still $1/2$ and the Hodge bound is $1/3$. Whether the infimum of $\F$ over maps falls below $1/2$ in the latter example remains open.

The four coordinates here give the smallest source dimension in which the harmonic representative can have rank greater than two. They are not a hypothesis of Theorem~\ref{thm:test}. Taking the product of the source with any connected closed Riemannian manifold gives examples in every dimension at least four. The lifted harmonic representative remains harmonic by the product Laplacian, and the energy and squared norms acquire the volume of the added factor. The horizontal vector fields in the proof of stationarity remain divergence-free, so the lifted projection is still a global symphonic minimizer. The small-oscillation test pulls back to a coclosed test of the same pointwise comass. Taking the target product with a Riemannian manifold and pulling back its torus area form likewise preserves the rank cap; projection onto the torus cannot increase the symphonic energy since the extra pullback metric is positive semidefinite.

\paragraph{Scope.}
Proposition~\ref{prop:quotient} gives a necessary and sufficient condition for strict improvement by the optimized test, not for a strict gap in every variational problem. The criterion uses the metric on the source and is not an invariant of the underlying topological manifold alone. On the same projection class, setting $b=0$ or taking $f$ constant makes the excess vanish. Determining whether pullbacks attain the convex relaxation outside the certificate regime is a separate question.

\paragraph{Use of computational assistance.}
OpenAI ChatGPT (Codex) was used to develop the rank estimate and shear example, draft and revise the proofs, search for related literature, and prepare the manuscript. The mathematical arguments above do not rely on numerical evidence.

\end{document}